\documentclass[12pt]{article}
\usepackage{latexsym}
\usepackage{amssymb}
\usepackage{graphicx}
\usepackage{cite}

\newtheorem{Theorem}{Theorem}[part]
\newtheorem{Definition}{Definition}[part]

\newtheorem{Lemma}{Lemma}[part]
\newtheorem{Corollary}{Corollary}[part]

\def \ep{\hbox{ }\hfill$\Box$}

\def \ra{\rightarrow}
\def\reff#1{{\rm(\ref{#1})}}

\begin{document}
\title{The E-Eigenvectors of Tensors}

\author{Shenglong Hu \thanks{Email: Tim.Hu@connect.polyu.hk. Department of
Applied Mathematics, The Hong Kong Polytechnic University, Hung Hom,
Kowloon, Hong Kong.}\hspace{4mm}\mbox{and}\hspace{4mm} Liqun Qi \thanks{Email:
maqilq@polyu.edu.hk. Department of Applied Mathematics, The Hong
Kong Polytechnic University, Hung Hom, Kowloon, Hong Kong. This
author's work was supported by the Hong Kong Research Grant
Council (Grant No. PolyU 501909, 502510, 502111 and 501212).}}

\date{\today}
\maketitle



\begin{abstract}
\noindent
We first show that the eigenvector of a tensor is well-defined. The differences between the eigenvectors of a tensor and its E-eigenvectors are the eigenvectors on the nonsingular projective variety $\mathbb S=\{\mathbf x\in\mathbb P^n\;|\;\sum\limits_{i=0}^nx_i^2=0\}$. We show that a generic tensor has no eigenvectors on $\mathbb S$. Actually, we show that a generic tensor has no eigenvectors on a proper nonsingular projective variety in $\mathbb P^n$. By these facts, we show that the coefficients of the E-characteristic polynomial are algebraically dependent. Actually, a certain power of the determinant of the tensor can be expressed through the coefficients besides the constant term.
Hence, a nonsingular tensor always has an E-eigenvector. When a tensor $\mathcal T$ is nonsingular and symmetric, its E-eigenvectors are exactly the singular points of a class of hypersurfaces defined by $\mathcal T$ and a parameter. We give explicit factorization of the discriminant of this class of hypersurfaces, which completes Cartwright and Strumfels' formula.  We show that the factorization contains the determinant and the E-characteristic polynomial of the tensor $\mathcal T$ as irreducible factors.

\vspace{3mm}

\noindent {\bf Key words:}\hspace{2mm} Tensor, E-eigenvector, E-characteristic polynomial, invariants, nonsingular
\vspace{3mm}

\noindent {\bf AMS subject classifications (2010):}\hspace{2mm} 14D21; 15A18; 15A69; 15A72
  \vspace{3mm}


\end{abstract}

\section{Introduction}
\setcounter{Theorem}{0} \setcounter{Proposition}{0}
\setcounter{Corollary}{0} \setcounter{Lemma}{0}
\setcounter{Definition}{0} \setcounter{Remark}{0}
\setcounter{Algorithm}{0}  \setcounter{Example}{0}
Let $\mathbb C$ be the field of complex numbers and $\mathbb C^{n+1}$ the ($n+1)$-dimensional complex space.
For a tensor $\mathcal T\in\otimes^m(\mathbb C^{n+1})$ of order $m$ and dimension $n+1$ with integers $m\geq 3,n\geq 1$, we identify it with its coefficient hypermatrix $(t_{i_1i_2\cdots i_m})$ under the canonical bases of $\otimes^m(\mathbb C^{n+1})$. Obviously, $t_{i_1\cdots i_m}\in\mathbb{C}$ for all $i_j\in(n):=\{0,1,\ldots,n\}$ and $j\in[m]:=\{1,\ldots,m\}$. Given a vector $\mathbf{x}\in \mathbb{C}^{n+1}$, define an ($n+1$)-dimensional vector ${\cal T}\mathbf{x}^{m-1}$ with its $i$-th element being $\sum\limits_{i_2,\ldots,i_m\in(n)}t_{ii_2\cdots i_m}x_{i_2}\cdots x_{i_m}$ for all $i\in(n)$.

\begin{Definition}\label{def-7}
Given $\mathcal T\in\otimes^m(\mathbb C^{n+1})$, a vector $\mathbf x\in\mathbb C^{n+1}\setminus\{0\}$ is called an eigenvector of the tensor $\mathcal T$ if
\begin{eqnarray}\label{eign}
\mathcal T\mathbf x^{m-1}\wedge \mathbf x=0.
\end{eqnarray}
\end{Definition}

The eigenvectors of tensors are proposed by Lim and Qi independently in \cite{l05,q}. There are several kinds of eigenvectors in the literature, see \cite{q,q2,l05,hhlq11,cs,nqww} and references therein. For the applications of eigenvectors, see \cite{hq12,oo12} and references therein. Definition \ref{def-7} comes from the E-eigenvectors of tensors introduced in \cite{q}. The precise definition of E-eigenvector is as follows: A number $\lambda\in\mathbb{C}$ is called an E-eigenvalue of ${\cal T}$, if it, together with a nonzero vector $\mathbf{x}\in\mathbb{C}^{n+1}$, satisfies
\begin{eqnarray}\label{e-eig}
\left\{\begin{array}{ccc}{\cal T}\mathbf{x}^{m-1}&=&\lambda \mathbf{x},\\\mathbf{x}^T\mathbf{x}&=&1.\end{array}\right.
\end{eqnarray}
$\mathbf{x}$ is then called the associated E-eigenvector of the E-eigenvalue $\lambda$, and $(\lambda,\mathbf{x})$ is called an eigenpair.
When $m=2$, this definition differs from the classic one for a matrix: it excludes the eigenpairs with eigenvectors $\mathbf x$ such that $\mathbf x^T\mathbf x=0$. Cartwright and Sturmfels \cite{cs} suggest the removing of the normalization $\mathbf x^T\mathbf x=1$ in \reff{e-eig} and use equivalence classes of eigenpairs. In the paper, we use Definition \ref{def-7}, which is equivalent to Cartwright and Sturmfels's equivalence classes of eigenpairs. This formulation is proposed by Oeding and Ottaviani in \cite{oo12}. With this definition, it is also unambiguous to say that a point $[\mathbf x]\in\mathbb P^n:=\mathbb P(\mathbb C^{n+1})$ is an eigenvector, which would mean that its representation $\mathbf x$ is an eigenvector in the sense of Definition \ref{def-7}.

Note that, by Definition \ref{def-7}, a nonzero vector $\mathbf x$ is an eigenvector of $\mathcal T$ if and only if $\mathcal T\mathbf x^{m-1}=\lambda \mathbf x$ for some $\lambda\in\mathbb C$. Since $\alpha \mathbf x$ is an eigenvector for any nonzero $\alpha\in\mathbb C$ whenever $\mathbf x$ is, it does not make sense to talk about eigenvalues when $\lambda\neq 0$. While, it does make sense to define the zero eigenvalue, since it keeps zero when its eigenvectors are scaled.
In \cite{cs}, if an eigenvector satisfying the normalization in \reff{e-eig}, then the corresponding $\lambda$ is called a normalized eigenvalue. It is easy to see that normalized eigenvalues are just E-eigenvalues.

In the investigation of eigenvalue theory of tensors, the theory of the characteristic polynomials of tensors plays a fundamental role, see \cite{q,lqz,hq13,hhlq11} and references therein. The E-characteristic polynomial introduced by Qi \cite{q} accompanies the E-eigenvalues of a tensor.

\begin{Definition}\label{def-2}
Given $\mathcal T\in\otimes^m(\mathbb C^{n+1})$, the E-characteristic polynomial $\chi_{{\cal T}}(\lambda)$ is defined as
\begin{eqnarray}\label{cha}
\chi_{{\cal T}}(\lambda):=\left\{\begin{array}{ll}\mbox{Res}_{\mathbf{x}}\left({\cal T}\mathbf{x}^{m-1}-\lambda(\mathbf{x}^T\mathbf{x})^{\frac{m-2}{2}}\mathbf{x}\right)&\mbox{when}\; m\;\mbox{is even},\\
\mbox{Res}_{\mathbf{x},\beta}\left(\begin{array}{c}{\cal T}\mathbf{x}^{m-1}-\lambda \beta^{m-2}\mathbf{x}\\\mathbf{x}^T\mathbf{x}-\beta^2\end{array}\right)&\mbox{when}\; m\;\mbox{is odd}.\end{array}\right.
\end{eqnarray}
Here $\mbox{Res}$ is the resultant of system of polynomials in the classic sense \cite{clo1,clo,gkz}.
\end{Definition}

The determinant of a tensor is an analogue of the determinant of a matrix, see \cite{hhlq11} for the definition and its various properties. By the determinant, we can define singular tensors.

\begin{Definition}\label{def-3}
Given $\mathcal T\in\otimes^m(\mathbb C^{n+1})$, it is called singular if $\mbox{Det}(\mathcal T)=0$; otherwise, it is called nonsingular.
\end{Definition}
This definition is used in \cite{h12,fo12}. It is different from the definitions by Qi \cite{q2}, Lim \cite{l05} and Cartwright and Sturmfels \cite{cs}. The latter three are different from each other as well. We see that Definition \ref{def-3} is an analogue of the definition for matrices. In the sequel, we can see that it is convenient to state results with this definition.

From \cite{hhlq11}, we get that $\mathcal T$ is singular if and only if there is a nonzero $\mathbf x\in\mathbb C^{n+1}$ such that $\mathcal T\mathbf x^{m-1}=0$. This fact will be used many times in the subsequent analysis.

Recently, the progress on the investigation for eigenvectors of tensors is great, see \cite{lqz,hq13,oo12,fo12,cs}. However, there are still some fundamental but unsolved problems. In this paper, we will concentrated on some of them as the followings:
\begin{itemize}
\item [1.] Is there always an eigenvector for an arbitrary tensor $\mathcal T\in\otimes^m(\mathbb C^{n+1})$?
\item [2.] Is there always an E-eigenvector for an arbitrary tensor $\mathcal T\in\otimes^m(\mathbb C^{n+1})$?
\item [3.] How many eigenvectors of a tensor that are not E-eigenvectors? And more generally, how many eigenvectors of a tensor are in a given nonsingular projective variety in $\mathbb P^n$?
\item [4.] There are $\frac{(m-1)^{n+1}-1}{m-2}$ E-eigenvalues for a generic tensor, while $\frac{(m-1)^{n+1}-1}{m-2}+1$ coefficients of the E-characteristic polynomial \cite{lqz,hq13,cs,q2}. Since both the E-eigenvalues and the coefficients are invariants of the underlying tensor, and they are related, is there algebraic dependence of these coefficients?
\end{itemize}
We show that these questions are related. We give affirmative answers to Questions 1 and 4, and answers for generic tensors to Questions 2 and 3. We remark that in \cite{q}, Questions 1 and 2 are solved with positive answers for real symmetric tensors of even orders. \cite[Corollar 2.3]{cs} solves Question 2 with a positive answer for real tensors when either $m$ is odd or $n$ is even. For the first half of Question 3, the proof of \cite[Lemma 7.2]{lqz} says that for a generic tensor of dimension two, it has at most two eigenvectors that are not E-eigenvectors.

The rest of this paper is organized as follows. In the next section, we first give some basic facts about the E-characteristic polynomials which are used extensively in the subsequent analysis. Then, we present the definition of vector bundles on a variety and the definition of sections.
In Section 3, based on a result from complex dynamics on projective spaces by Fornaess and Sibony \cite{fs}, we show that every given tensor has at least one eigenvector, i.e., the eigenvector is well-defined. Then, in Section 4, we discuss the eigenvectors on nonsingular projective varieties with the help of the concept of vector bundles and a recent result by Friedland and Ottaviani \cite{fo12}. We prove that a generic tensor has no eigenvectors on the variety $\mathbb S=\{\mathbf x\in\mathbb P^n\;|\;\sum\limits_{i=0}^nx_i^2=0\}$. With these results,
we prove that the coefficients of the E-characteristic polynomial are algebraically dependent in Section 5. As a consequence, a nonsingular tensor always has an E-eigenvector. In Section 6, we then investigate the E-eigenvectors of a nonsingular symmetric tensor. We complement a formula for the discriminant, initiated by Cartwright and Sturmfels, of the affine hypersurfaces defined by a symmetric tensor and a parameter. We show that the discriminant of this class of hypersurfaces contains the determinant and the E-characteristic polynomial of the tensor as irreducible factors. We conclude this paper in Section 7 with some final remarks and extensions on the expression of the determinant as a combination of invariants of the underlying tensor.

\section{Preliminaries}
\setcounter{Theorem}{0} \setcounter{Proposition}{0}
\setcounter{Corollary}{0} \setcounter{Lemma}{0}
\setcounter{Definition}{0} \setcounter{Remark}{0}
\setcounter{Algorithm}{0}  \setcounter{Example}{0}
\hspace{4mm}
In this section, some basic facts about the E-characteristic polynomials and eigenvectors are presented. We also give the basic definitions of vector bundles and sections on varieties.

\subsection{The E-Characteristic Polynomial}
We first present some notation used in this paper, which would be clear from the introduction already.
Scalars are written as lowercase letters ($\lambda,a,\ldots$), vectors are written as bold lowercase letters ($\mathbf{x},\mathbf{y},\ldots$), the $i$-th entry of a vector $\mathbf{x}$ is denoted by $x_i$, matrices are written as italic capitals ($A,B,\ldots$), tensors are written as calligraphic capitals (${\cal T}$, ${\cal D}$, $\ldots$), and sets are written as black bold letters ($\mathbb M$, $\mathbb S$, $\ldots$).

Given a ring $\mathbb{K}$ (hereafter, we mean a commutative ring with $1$),
we denote by $\mathbb{K}[\mathbb{E}]$ the polynomial ring consisting of polynomials in the set $\mathbb{E}$ of indeterminate variables with coefficients in $\mathbb{K}$. Especially, we denote by $\mathbb{K}[{\cal T}]$ the polynomial ring consisting of polynomials in indeterminate variables $\{t_{i_1\ldots i_m}\}$ with coefficients in $\mathbb{K}$, and similarly for $\mathbb{K}[\lambda]$, $\mathbb{K}[\lambda,{\cal T}]$, etc.

We summarize some important results on eigenvectors and the E-characteristic polynomials of tensors in the following theorem.

\begin{Theorem}\label{thm-1}
Given $\mathcal T\in\otimes^m(\mathbb C^{n+1})$ and its E-characteristic polynomial $\chi_{\mathcal T}(\lambda)$, then, we have the followings.
\begin{itemize}
\item [(i)] The tensor $\mathcal T$ has zero as its eigenvalue if and only if $\mathcal T$ is singular.
\item [(ii)] If $\mathcal T$ is nonsingular, then $\lambda\in\mathbb C$ is an E-eigenvalue of the tensor $\mathcal T$ if and only if $\lambda$ is a root of $\chi_{\mathcal T}(\lambda)$.
\item [(iii)] Let $\mathbb A^1$ be the set of tensors having all complex numbers as their E-eigenvalues; $\mathbb A^2$ be the set of tensors having infinitely many E-eigenvalues; $\mathbb A^3$ the set of tensors having zero E-characteristic polynomials; and $\mathbb A^4$ the set of singular tensors. Then, we have the strict inclusions
    \begin{eqnarray}\label{set}
    \mathbb A^4\subset\mathbb A^3\subset \mathbb A^2\subset \mathbb A^1.
    \end{eqnarray}
\item [(iv)] The degree $\mbox{deg}(\psi_{{\cal T}}(\lambda))$ is equal to $\frac{(m-1)^{n+1}-1}{m-2}$ for generic even order tensors and it is equal to $2\frac{(m-1)^{n+1}-1}{m-2}$ for generic odd order tensors.
\item [(v)] If the tensor $\mathcal T$ has only finitely many eigenvectors (up to scalar multiplication), then the number of them, counted with multiplicity, is equal to $\frac{(m-1)^{n+1}-1}{m-2}$.
\item [(vi)] If the tensor $\mathcal T$ is symmetric, then the number of E-eigenvalues of $\mathcal T$ is finite.
\item [(vii)] The coefficients of the E-characteristic polynomial of a tensor are invariants under the action of the group of orthogonal matrices in the sense of \cite{q,lqz}.
\item [(viii)] The constant term of the E-characteristic polynomial of the tensor $\mathcal T$ is equal to $\mbox{Det}(\mathcal T)$ when $m$ is even and $\mbox{Det}(\mathcal T)^2$ when $m$ is odd.
\end{itemize}
\end{Theorem}

\noindent {\bf Proof.} (i) follows from \cite[Theorem 3.1(i)]{hhlq11}.

(ii) If $(\lambda, \mathbf x)$ is an eigenpair of $\mathcal T$, then $\lambda$ would be a root of $\chi_{\mathcal T}(\lambda)$, since the corresponding polynomial system in \reff{cha} has a nontrivial solution. On the contrary, if the polynomial system in \reff{cha} has a nontrivial solution $\mathbf x$, then $\mathbf x^T\mathbf x\neq 0$. Since otherwise we have $\mathcal T\mathbf x^{m-1}=0$, which is a contradiction to the nonsingularity. Consequently, we can normalize $\mathbf x$ and then it is easy to see that $\lambda$ is an E-eigenvalue.

(iii) follows from \cite[Proposition 3.4, Examples 3.2 and 3.5]{cs} and \cite[Eq. (9)]{q2}. An explicit singular tensor with nonzero E-characteristic polynomial can be easily found by the formulae in \cite[Example 3.5]{cs}.

(iv) follows from Definition \ref{def-2}, and \cite[Theorem 1.2]{cs} which also implies (v).

(vi) follows from \cite[Theorem 5.6]{cs}.

(vii) follows from \cite[Theorem 3.3]{lqz}.

Finally, (viii) follows from \cite[Theorem 7.4]{lqz} and \cite[Theorems 3.1 and 4.1]{hq13}. The proof is complete. \ep

The following lemma is \cite[Proposition 2.4(ii)]{hhlq11}.
\begin{Lemma}\label{lem-0}
If we view $\mathcal T$ as an $m$-th order ($n+1$)-dimensional tensor consisting of indeterminate variables $t_{i_1\ldots i_m}$, then $\mbox{Det}(\mathcal T)$ is an irreducible homogeneous polynomial in $\mathbb C[\mathcal T]$ of degree $(n+1)(m-1)^n$.
\end{Lemma}

We also have the following lemma.
\begin{Lemma}\label{lem-00}
If we view $\mathcal T$ as an $m$-th order ($n+1$)-dimensional tensor consisting of indeterminate variables $t_{i_1\ldots i_m}$ and $\lambda$ an indeterminate variable, then $\chi_{\mathcal T}(\lambda)$ is an irreducible homogeneous polynomial in $\mathbb C[\mathcal T,\lambda]$ of degree $(n+1)(m-1)^n$. Moreover, for any generic realization $\mathcal T$ in $\otimes^m(\mathbb C^{n+1})$, $\chi_{\mathcal T}(\lambda)\in\mathbb C [\lambda]$ is irreducible.
\end{Lemma}

\noindent {\bf Proof.} Suppose on the contrary that $\chi_{\mathcal T}(\lambda)$ can be reduced as $\chi_{\mathcal T}(\lambda)=p(\mathcal T,\lambda)q(\mathcal T,\lambda)$. We note that both $p$ and $q$ should be homogeneous as polynomials in $\mathbb C[\mathcal T,\lambda]$. If the degree of both $p$ and $q$ with respect to $\lambda$ is larger than zero, then $\chi_{\mathcal T}(\lambda)\in(\mathbb C[\mathcal T])[\lambda]$ can be reduced, which contradicts \cite[Corollary 3.1]{cs}. Then, the only chance for a reduced factorization could be that one of $p$ and $q$ is independent of $\lambda$. Suppose $p(\mathcal T,\lambda)\in\mathbb C[\mathcal T]$, without loss of generality. Then, it is safe to write $p(\mathcal T)$ for $p(\mathcal T,\lambda)$.

When $m$ is even, we must have $p(\mathcal T)=1$ by Lemma \ref{lem-0} and Theorem \ref{thm-1} (viii). Since otherwise, $\mbox{Det}(\mathcal T)$ would be reducible. By the same reason, $p(\mathcal T)=\mbox{Det}(\mathcal T)$ when $m$ is odd. Consequently, $\chi_{\mathcal T}(\lambda)=\mbox{Det}(\mathcal T)q(\mathcal T,\lambda)$. This contradicts the last strict inclusion in \reff{set} of Theorem \ref{thm-1} (iii).

The last conclusion follows from \cite[Corollary 3.1]{cs}. The proof is complete. \ep

\subsection{Vector Bundles}
In this subsection, we give the definitions of vector bundles. For comprehensive references and the definition of varities, please see \cite{gh,h77,oss80,skkt00}.

\begin{Definition}\label{def-5}
A vector bundle of rank $r$ on a variety $\mathbb{X}$ is a variety $\mathbb{E}$, together with a morphism $\pi: \mathbb{E}\ra\mathbb{X}$, which is called the projection, such that the following conditions are satisfied:
\begin{itemize}
\item [(i)] There is an open cover $\cup\mathbb{U}_i$ of $\mathbb{X}$ such that $\pi^{-1}(\mathbb{U}_i)$ is isomorphic to the product $\mathbb{U}_i\times\mathbb{C}^r$ by fiber-preserving maps, i.e., there are isomorphisms $\varphi_i: \pi^{-1}(\mathbb{U}_i)\ra \mathbb{U}_i\times\mathbb{C}^r$, which are called trivializations, such that $\pi=p\circ\varphi_i$ on $\pi^{-1}(\mathbb{U}_i)$. Here $p: \mathbb{U}_i\times\mathbb{C}^r\ra \mathbb{U}_i$ is the natural projection onto the first factor.
\item [(ii)] The isomorphisms $\varphi_i$ are linearly compatible in the following sense: On $\mathbb{U}_i\cap \mathbb{U}_j$, the composition
\begin{eqnarray*}
\varphi_j\circ\varphi_i^{-1}: (\mathbb{U}_i\cap\mathbb{U}_j)\times\mathbb{C}^r&\ra& (\mathbb{U}_i\cap\mathbb{U}_j)\times\mathbb{C}^r\\
(\mathbf{x},\mathbf{v})&\mapsto& (\mathbf{x}, q\circ \varphi_j\circ\varphi_i^{-1}\circ (\mathbf{x},\mathbf{v})),
\end{eqnarray*}
is a linear map with respect to $\mathbb{C}^r$ for any fixed value of $\mathbf{x}$. Here $q: \mathbb{U}_i\times\mathbb{C}^r\ra \mathbb{C}^r$ is the natural projection onto the second factor.
\end{itemize}
\end{Definition}

From the definition, we see that $\pi^{-1}(\mathbf x)$ is isomorphic to a vector space of dimension $r$. We denoted this vector space as $\mathbb E_{\mathbf x}$.
The variety $\mathbb{E}$ is called the total space of the vector bundle. Sometimes, the entire vector bundle is denoted simply by its total space. Vector bundle of rank one is called line bundle.

\begin{Definition}\label{def-6}
Let $\pi: \mathbb{E}\ra\mathbb{X}$ be a vector bundle, and $\mathbb{U}\subseteq\mathbb{X}$ be an open set. A section of the vector bundle on the set $\mathbb{U}$ is a morphism $s: \mathbb{U}\ra\mathbb{E}$ such that $\pi\circ s$ is the identity map on $\mathbb{U}$. The set of all sections of $\mathbb{E}$ over $\mathbb{U}$ is denoted by $\mathcal{E}(\mathbb{U})$. The set of global sections is defined to be the set $\mathcal{E}(\mathbb{X})$ of sections of $\mathbb{E}$ over the whole variety $\mathbb{X}$. It is denoted by $H^0(\mathbb E, \mathbb X)$. The zero locus of a section $s$ is the set of points $\mathbf x$ for which $s(\mathbf x)$ is the zero vector in $\mathbb E_{\mathbf x}$.
\end{Definition}

A subspace $\mathbb V\subset H^0(\mathbb E, \mathbb X)$ is said to generate $\mathbb E$ if the set $\{v(\mathbf x)\;|\;v\in\mathbb V\}$ is equal to $\pi^{-1}(\mathbf x)$ at every $\mathbf x\in\mathbb X$.

In this paper, we mainly consider vector bundles on projective varieties. For example, $\mathbb C^{n+1}$ can be viewed as a vector bundle on the projective space $\mathbb P^n$ in a natural way. Take $\pi^{-1}:\mathbb P^n\ra \mathbb C^{n+1}$ as $\pi^{-1}([\mathbf x])=\langle\mathbf x\rangle$. This is a line bundle, and is usually called the tautological line bundle. We denoted it by $\mathcal O(-1)$.

By Chow's theorem \cite{c49,h77,skkt00}, we have that every compact complex manifold embedded in $\mathbb P^n$ is a projective variety. Consequently, we can define vector bundles on compact complex manifolds in $\mathbb P^n$ with the above definitions.

\section{Well-Definiteness of Eigenvectors}
\setcounter{Theorem}{0} \setcounter{Proposition}{0}
\setcounter{Corollary}{0} \setcounter{Lemma}{0}
\setcounter{Definition}{0} \setcounter{Remark}{0}
\setcounter{Algorithm}{0}  \setcounter{Example}{0}
\hspace{4mm}
In this section, we discuss the well-definiteness of the eigenvectors defined by Definition \ref{def-7}.
By well-definiteness of the eigenvector definition (Definition \ref{def-7}), we mean that there exists an eigenvector for any given tensor.

To this end, we give some equivalent reformulations of the eigenvectors of tensors first.

\begin{Theorem}\label{thm-2}
Given $\mathcal T\in\otimes^m(\mathbb C^{n+1})$, then under the equivalence of nonzero scalar multiplication, the following three sets are the same.
\begin{itemize}
\item [(1)] The set of eigenvectors of the tensor $\mathcal T$.
\item [(2)] The union of the sets of nonzero solutions of the following two polynomial systems:
\begin{eqnarray}
\mathcal T\mathbf x^{m-1}=\mathbf x\label{unit},\\
\mathcal T\mathbf x^{m-1}=\mathbf 0\label{zero}.
\end{eqnarray}
\item [(3)] The set of nonzero solutions of the following polynomial system:
\begin{eqnarray}\label{homo}
x_i\left(\mathcal T\mathbf x^{m-1}\right)_j-x_j\left(\mathcal T\mathbf x^{m-1}\right)_i=0,\;0\leq i<j\leq n.
\end{eqnarray}
\end{itemize}
\end{Theorem}

\noindent {\bf Proof.} All the equivalences follow from the eigenvector definition (Definition \ref{def-7}). For example, \reff{homo} follows from the elimination of the variable $\lambda$ from the equivalent reformulation as $\mathcal T\mathbf x^{m-1}=\lambda\mathbf {x}$.  \ep

The equivalent reformulations given above provide new ways to investigate the eigenvectors of tensors. From \reff{homo}, it is clear that the set of eigenvectors is a projective set. \reff{zero} says that the set of eigenvectors corresponding to the zero eigenvalue of a symmetric tensor is the set of singular points of a projective hypersurface. From \reff{unit}, we see that the eigenvectors of a nonsingular tensor are just the fixed points of a dynamic system defined on the projective space. Actually, we will prove the well-definiteness of Definition \ref{def-7} by using a result from dynamics on complex projective spaces, which is due to Fornaess and Sibony \cite{fs}.

Let $\mathbb P^n:=\mathbb P(\mathbb C^{n+1})$ be the complex projective space of $\mathbb C^{n+1}$. Then, given $\mathcal T\in\otimes^m(\mathbb C^{n+1})$, we can define a dynamic $f: \mathbb{P}^{n}\ra\mathbb{P}^{n}$ as
\begin{eqnarray*}
f([\mathbf x]):=[\mathcal T\mathbf x^{m-1}].
\end{eqnarray*}
We see that $f$ is defined everywhere on $\mathbb P^n$ if and only if $\mathcal T$ is nonsingular. Actually, we have the following result.
\begin{Lemma}\label{lem-000}
Let $\mathcal T\in\otimes^m(\mathbb C^{n+1})$ be nonsingular and $f$ be defined as above, then $f$ is holomorphic and it has $\frac{(m-1)^{n+1}-1}{m-2}$ fixed points (with multiplicity).
\end{Lemma}

\noindent {\bf Proof.} By the definition of $f$, it is easy to see that it is holomorphic. Consequently, by \cite[Corollary 3.2]{fs}, the result on the number of fixed points follows. \ep

By Theorem \ref{thm-2} and Lemma \ref{lem-000}, we can show that the eigenvector definition is well-defined.
\begin{Theorem}\label{thm-3}
The eigenvector definition (Definition \ref{def-7}) is well-defined.
\end{Theorem}

\noindent {\bf Proof.} If the tensor $\mathcal T$ is singular, then by Theorem \ref{thm-1}(i) the system \reff{zero} has a nontrivial solution. Consequently, by Theorem \ref{thm-2}, we know that the eigenvector definition (Definition \ref{def-7}) is well-defined.

If the tensor $\mathcal T$ is nonsingular, by Lemma \ref{lem-000}, $f$ has $\frac{(m-1)^{n+1}-1}{m-2}$ fixed points (counted with multiplicity). Note that the fixed points of the map $f$ are exactly the nontrivial solutions of the system \reff{unit} for the tensor $\mathcal T$. Consequently, the eigenvector definition (Definition \ref{def-7}) is well-defined by Theorem \ref{thm-2}. \ep

This gives an affirmative answer to Question 1.

\section{Eigenvectors of Tensors on Nonsingular Projective Varieties}
\setcounter{Theorem}{0} \setcounter{Proposition}{0}
\setcounter{Corollary}{0} \setcounter{Lemma}{0}
\setcounter{Definition}{0} \setcounter{Remark}{0}
\setcounter{Algorithm}{0}  \setcounter{Example}{0}
\hspace{4mm}
Let $\mathbb C^{n+1}$ be the ($n+1$)-dimensional complex space and $\mathbb P^n$ its projective space. Let $\mathbb M$ be a nonsingular projective variety in $\mathbb P^n$. A nonsingular variety is a variety whose tangent space at each point has the same dimension. We denote by $\mbox{dim}(\mathbb M)$ the dimension of $\mathbb M$. In this section, we investigate eigenvectors of a given tensor $\mathcal T\in\otimes^m(\mathbb C^{n+1})$ in a given nonsingular projective variety $\mathbb M\subseteq\mathbb P^n$. Especially, we will show that a generic tensor has no eigenvectors on $\mathbb S:=\{[\mathbf x]\in\mathbb P^n\;|\;\sum\limits_{i\in(n)}x_i^2=0\}$. Here $(x_0:x_1:\ldots:x_n)$ are the homogeneous coordinates of the point $[\mathbf x]\in\mathbb P^n$. We refer to \cite{clo,clo1,f84,gh,h77} for the basic notation and concepts.

Denote by $\mathcal O(-1)$ the tautological line bundle on $\mathbb P^n$ and $\mathcal O(1)$ its dual. Denote by $\mathbb Q$ the quotient bundle $(\mathcal O\otimes \mathbb C^{n+1})/\mathcal O(-1)$. Then, we have an exact sequence
\begin{eqnarray*}
0\ra\mathcal O(-1)\ra \mathcal O\otimes \mathbb C^{n+1}\ra \mathbb Q\ra 0.
\end{eqnarray*}
From this sequence, we can derive the number of eigenvectors of tensors by Chern's classes, see \cite{oo12,fo12}. This method possesses great success recently \cite{os12,oo12,fo12}.

Similarly, we can defined vector bundles $\mathcal O_{\mathbb M}(-1)$ and $\mathbb Q_{\mathbb M}:=(\mathcal O_{\mathbb M}\otimes \mathbb C^{n+1})/\mathcal O_{\mathbb M}(-1)$ on the nonsingular projective variety $\mathbb M$.
We also have an exact sequence
\begin{eqnarray*}
0\ra\mathcal O_M(-1)\ra \mathcal O_M\otimes \mathbb C^{n+1}\ra \mathbb Q_{\mathbb M}\ra 0.
\end{eqnarray*}

For $[\mathbf x]\otimes \mathbf y\in\mathcal O\otimes\mathbb C^{n+1}$, denoted by $[[\mathbf y]]\in \mathbb Q_{[\mathbf x]}$ the image of the projection onto $\mathbb Q_{[\mathbf x]}$ of the natural morphism $\mathcal O\otimes\mathbb C^{n+1}\ra \mathbb Q$.
Denote by $\mathbb E$ the vector bundle $\mathbb Q_{\mathbb M}$ on $\mathbb M$ for the convenience. Then, $\mbox{rank}(\mathbb E)=n$ by Definition \ref{def-5}.
Define a map $L: \otimes^m(\mathbb C^{n+1})\ra H^0(\mathbb E,\mathbb M)$ as
\begin{eqnarray}\label{sect}
L(\mathcal T)([\mathbf x]):=[\mathbf x]\otimes[[\mathcal T\mathbf x^{m-1}]],
\end{eqnarray}
where $[[\mathcal T\mathbf x^{m-1}]]\in (\mathbb Q_{\mathbb M})_{[\mathbf x]}$. It is easy to see that $L(\mathcal T)$ is a global section of $\mathbb E$ by Definition \ref{def-6}. The next lemma shows that $L(\otimes^m(\mathbb C^{n+1}))$ generates $\mathbb E$.

\begin{Lemma}\label{lem-1}
Let $L$ be defined as above and
let $\mathbb V:=\{L(\mathcal T)\;|\;\mathcal T\in\otimes^m(\mathbb C^{n+1})\}$. Then, $\mathbb V$ is a subspace of the linear space $H^0(\mathbb E,\mathbb M)$ of all global sections on $\mathbb E$. Moreover, $\mathbb V$ generates $\mathbb E$.
\end{Lemma}

\noindent {\bf Proof.}
Note that $[[(\mathcal T+\mathcal U)\mathbf x^{m-1}]]=[[\mathcal T\mathbf x^{m-1}]]+[[\mathcal U\mathbf x^{m-1}]]$ by the definition of the quotient space $(\mathbb Q_{\mathbb M})_{[\mathbf x]}$. Consequently, it is easy to see that $\mathbb V$ is a linear subspace of $H^0(\mathbb E, \mathbb M)$.

In the following, for any $[[\mathbf y]]\in (\mathbb Q_{\mathbb M})_{[\mathbf x]}$,  defined a tensor $\mathcal T$ as
\begin{eqnarray*}
\frac{1}{\alpha^{m-1}}\mathbf y\otimes(\otimes^{m-1}\mathbf x^*),
\end{eqnarray*}
where $\mathbf x^*$ is the conjugate of $\mathbf x$.
Denote $\alpha:=\mathbf x^T\mathbf x^*\neq 0$. We have that $L(\mathcal T)([\mathbf x])=[\mathbf x]\otimes[[\alpha^{m-1}\frac{1}{\alpha^{m-1}}\mathbf y]]=[\mathbf x]\otimes[[\mathbf y]]$. Consequently, $\mathbb V$ generates $\mathbb E$.

The proof is complete. \ep

We have a simple lemma as follows.
\begin{Lemma}\label{lem-0000}
Let $\mathbb M$ be a nonsingular projective variety in $\mathbb P^n$ with dimension $k$ and $\mathcal T\in\otimes^m(\mathbb C^{n+1})$. Let the map $L$ be defined as above. Then, $[\mathbf x]\in\mathbb M$ is an eigenvector of $\mathcal T$ if and only if it belongs to the zero locus of $L(\mathcal T)$.
\end{Lemma}

\noindent {\bf Proof.} If $\mathcal T \mathbf x^{m-1}\wedge \mathbf x=0$, then $[[\mathcal T \mathbf x^{m-1}]]=[[\lambda\mathbf x]]=0$ for some $\lambda\in\mathbb C$. Hence, $L(\mathcal T)([\mathbf x])=[\mathbf x]\otimes 0$ by the definition.

Conversely, if $L(\mathcal T)([\mathbf x])=0$, then $[[\mathcal T \mathbf x^{m-1}]]=0$. Since the kernel of the natural morphism $\mathcal O_M\otimes\mathbb C^{n+1}\ra \mathbb Q_{\mathbb M}$ is $\mathcal O_M(-1)$, we conclude that $\mathcal T \mathbf x^{m-1}\wedge \mathbf x=0$. Hence, $[\mathbf x]$ is an eigenvector of $\mathcal T$ by Definition \ref{def-7}. The proof is complete. \ep

Since $\mathbb P^n$ is a compact complex manifold of dimension $n$, every nonsingular projective variety $\mathbb X\subseteq \mathbb P^n$ can be regarded as a compact complex manifold of dimension $\mbox{dim}(\mathbb X)$. The coordinate chart can be obtained through the implicit function theorem.
Thus, the following theorem is a direct consequence of \cite[Theorem 2]{fo12}, which is due to Friedland and Ottaviani.

\begin{Theorem}\label{thm-5}
Let $\mathbb E$ be a vector bundle on a nonsingular projective variety $\mathbb M$ and $\mathbb V\subset H^0(\mathbb E, \mathbb M)$  be a subspace which generates $\mathbb E$. If $\mbox{rank}(\mathbb E)>\mbox{dim}(\mathbb M)$, then for a generic $s\in \mathbb V$, the zero locus of $s$ is empty.
\end{Theorem}

The cases when $\mbox{rank}(\mathbb E)\leq\mbox{dim}(\mathbb M)$ are discussed in \cite{fo12,oo12} as well. For more details on this result, please check the generic smoothness theorem and Bertini's theorem \cite[Corollary III.10.7 and Theorem II. 8.1.8]{h77}.

By Lemmas \ref{lem-1} and \ref{lem-0000}, and Theorem \ref{thm-5}, we immediately have the following theorem.
\begin{Theorem}\label{thm-7}
Let $\mathbb M$ be a nonsingular projective variety in $\mathbb P^n$ with dimension $k<n$. Then, a generic $\mathcal T\in\otimes^m(\mathbb C^{n+1})$ has no eigenvectors on $\mathbb M$.
\end{Theorem}

\noindent {\bf Proof.} Let $\mathbb V$ be defined as that in Lemma \ref{lem-1}. By Lemmas \ref{lem-1} and \ref{lem-0000}, and Theorem \ref{thm-5}, it is sufficient to prove that there is an open set (in the Zariski sense) of $\mathbb V$ corresponding to an open set of $\otimes^m(\mathbb C^{n+1})$. Since the intersection of two open sets is still open, if we can prove that the map $L$ defined as \reff{sect} is injective on an open set of $\mathcal T\in\otimes^m(\mathbb C^{n+1})$, then the statement follows immediately.

If $k=0$, then by Theorem \ref{thm-1} (ii) and Lemma \ref{lem-00} two generic tensors has no common eigenvectors. Hence, the eigenvectors of generic tensors cannot be a fixed finite set.
If $k>1$, then by the proof of Theorem \ref{thm-3} a generic tensor has only finitely many eigenvectors. Hence, the eigenvectors of a generic tensor is a finite set. Consequently, for a generic tensor $\mathcal T$, there exists $[\mathbf x]\in\mathbb M$ such that $[[\mathcal T\mathbf x^{m-1}]]\neq 0$. Hence, we have that $L$ is injective for generic tensors. \ep

Since the variety $\{[\mathbf x]\in\mathbb P^n\;|\; \sum\limits_{i=1}^{n+1}x_i^2=0\}\subset\mathbb P^n$ is nonsingular and has dimension $n-1$, we get the following corollary as a direct consequence of Theorem \ref{thm-7}.

\begin{Corollary}\label{cor-1}
For a generic $\mathcal T\in\otimes^m(\mathbb C^{n+1})$, it has no eigenvectors on $\{[\mathbf x]\in\mathbb P^n\;|\; \sum\limits_{i=1}^{n+1}x_i^2=0\}$.
\end{Corollary}

Since E-eigenvectors corresponding to eigenvectors outside the set $\{[\mathbf x]\in\mathbb P^n\;|\; \sum\limits_{i=1}^{n+1}x_i^2=0\}$, this gives an answer to the first half of Question 3 for generic tensors. Theorem \ref{thm-7} gives an answer to the second half of Question 3 for generic tensors.

\section{Existence of E-Eigenvectors}
\setcounter{Theorem}{0} \setcounter{Proposition}{0}
\setcounter{Corollary}{0} \setcounter{Lemma}{0}
\setcounter{Definition}{0} \setcounter{Remark}{0}
\setcounter{Algorithm}{0}  \setcounter{Example}{0}
\hspace{4mm}
In this section, we discuss the existence of E-eigenvectors of a given tensor $\mathcal T\in\otimes^m(\mathbb C^{n+1})$. To this end, we investigate the algebraic dependence of the coefficients of the E-characteristic polynomial of this tensor. By the algebraic dependence, we show that a nonsingular tensor always has an E-eigenvector.

In the following, we view $\mathcal T$ as an $m$-th order ($n+1$)-dimensional tensor consisting of indeterminate variables $t_{i_1\ldots i_m}$. By Theorem \ref{thm-1} (iv), we let $a_i$ be the codegree $i$ coefficient of the E-characteristic polynomial of the tensor $\mathcal T$ for $i\in(\frac{(m-1)^{n+1}-1}{m-2})$ when $m$ is even, and the codegree $2i$ coefficient of the E-characteristic polynomial of the tensor when $m$ is odd.
Then, by Theorem \ref{thm-1} (viii), we have
\begin{eqnarray*}
\chi_{\mathcal T}(\lambda)=a_0\lambda^{\frac{(m-1)^{n+1}-1}{m-2}}+a_1\lambda^{\frac{(m-1)^{n+1}-1}{m-2}-1}+\ldots+a_{\frac{(m-1)^{n+1}-1}{m-2}-1}\lambda+\mbox{Det}(\mathcal T)
\end{eqnarray*}
when $m$ is even; and
\begin{eqnarray*}
\chi_{\mathcal T}(\lambda)=a_0\lambda^{2\frac{(m-1)^{n+1}-1}{m-2}}+a_1\lambda^{2\frac{(m-1)^{n+1}-1}{m-2}-2}+\ldots+a_{\frac{(m-1)^{n+1}-1}{m-2}-1}\lambda^2+[\mbox{Det}(\mathcal T)]^2
\end{eqnarray*}
when $m$ is odd. The expression for the odd case with only even powers of $\lambda$ follows from the discussion on \cite[Page 1371]{q2}.

By Lemma \ref{lem-00}, $a_i\in\mathbb C[\mathcal T]$ is homogeneous of degree $(n+1)(m-1)^n-\frac{(m-1)^{n+1}-1}{m-2}+i$ when $m$ is even, and $2\left[(n+1)(m-1)^n-\frac{(m-1)^{n+1}-1}{m-2}+i\right]$ when $m$ is odd. For an ideal $I$ in a polynomial ring, we denoted by $\mathbb V(I)$ the algebraic set determined by $I$ in the ambient space. Correspondingly, we denoted by $\mathbb I(S)$ the radical ideal determined by a set $S$ \cite{clo,clo1}.
As polynomials in the variable $\mathcal T$, we have the following theorem.

\begin{Theorem}\label{thm-4}
Let $a_i$ be defined as above. We have
\begin{eqnarray}\label{ideal}
\mbox{Det}(\mathcal T)\in \sqrt{\langle a_0,\ldots,a_{t}\rangle},
\end{eqnarray}
where $t=\frac{(m-1)^{n+1}-1}{m-2}-1$.
\end{Theorem}

\noindent {\bf Proof.} Let $V_1:=\mathbb V(\langle \mbox{Det}(\mathcal T)\rangle)$ and $V_2:=\mathbb V(\langle a_0,\ldots,a_{t}\rangle)$. We note that the ambient space is $\otimes^m(\mathbb C^{n+1})$ and the corresponding polynomial ring is $\mathbb C[\mathcal T]$. By Lemma \ref{lem-0} and Hilbert's Nullstellensatz (see \cite[Theorem 4.2]{clo}), it is sufficient to prove that
\begin{eqnarray*}
V_2\subseteq V_1.
\end{eqnarray*}
Actually, this implies that
\begin{eqnarray*}
\mbox{Det}(\mathcal T)\in \langle \mbox{Det}(\mathcal T)\rangle=\sqrt{\mathbb I(V_1)}=\mathbb I(V_1)\subseteq \mathbb I(V_2)=\sqrt{\langle a_0,\ldots,a_{t}\rangle}.
\end{eqnarray*}
Consequently, \reff{ideal} follows immediately.

By Corollary \ref{cor-1}, there exists an open set $U$ (in the Zariski sense) such that every tensor in $U$ has no eigenvectors on the projective variety $\{[\mathbf x]\in\mathbb{P}^{n}\;|\;\sum\limits_{i\in(n)} x_i^2=0\}$.

Let $\mathcal T\in V_2$. If $\mathcal T\notin V_1$, then $\mbox{Det}(\mathcal T)\neq 0$. By Theorem \ref{thm-1} (ii), $\lambda\in\mathbb C$ is an E-eigenvalue if and only if it is a root of the E-characteristic polynomial.
Since $\mathcal T\in V_2\setminus V_1$, all the coefficients but the constant term of the E-characteristic polynomial are zero. Consequently, the E-characteristic polynomial has no root. So, every eigenvector of the tensor $\mathcal T$ is in the projective variety $\{[\mathbf x]\in\mathbb{P}^{n}\;|\;\sum\limits_{i\in(n)} x_i^2=0\}$. As a consequence of Corollary \ref{cor-1}, such a tensor $\mathcal T$ must be in the complement of $U$. In other words, we must have
\begin{eqnarray*}
U\cap V_2\subseteq U\cap V_1.
\end{eqnarray*}
Taking the Zariski closure on both sides in the ambient space, we get that $V_2\subseteq V_1$. Actually, we have
\begin{eqnarray*}
\mathbb I(U\cap V_2)=\sqrt{\mathbb I(U)+\mathbb I(V_2)}=\mathbb I(V_2).
\end{eqnarray*}
Here the second equality follows from the fact that $U$ is open, and hence $\mathbb I(U)=0$. Similarly, we have $\mathbb I(U\cap V_1)=\mathbb I(V_1)$. Consequently, we have $\mathbb I(V_1)\subseteq \mathbb I(V_2)$ since $\mathbb I(U\cap V_1)\subseteq \mathbb I(U\cap V_2)$ by $U\cap V_2\subseteq U\cap V_1$. From this, we have
\begin{eqnarray*}
\mathbb V(\mathbb I(V_2))\subseteq \mathbb V(\mathbb I(V_1)).
\end{eqnarray*}
Since both $V_1$ and $V_2$ are varieties, we have $V_2=\mathbb V(\mathbb I(V_2))\subseteq \mathbb V(\mathbb I(V_1))=V_1$ as claimed.

The proof is complete. \ep

This theorem gives an affirmative answer to Question 4.
By Theorem \ref{thm-4}, we have the following result on the existence of E-eigenvectors of a nonsingular tensor.

\begin{Theorem}\label{thm-8}
If the tensor $\mathcal T$ is nonsingular, then it always has an E-eigenvector.
\end{Theorem}

\noindent {\bf Proof.} Let $a_i$ be defined as above, by \reff{ideal}, there exist positive integer $k$ and polynomials $b_i\in\mathbb C[\mathcal T]$ such that
\begin{eqnarray*}
[\mbox{Det}(\mathcal T)]^k=\sum\limits_{i\in(\frac{(m-1)^{n+1}-1}{m-2}-1)}a_ib_i.
\end{eqnarray*}
Consequently, there exists at least one $a_i\neq 0$ whenever $\mbox{Det}(\mathcal T)\neq 0$. By the fundamental theorem of algebra, we see that the E-characteristic polynomial $\chi_{\mathcal T}(\lambda)$ must have a root. By Theorem \ref{thm-1} (ii), this root is an E-eigenvalue of the tensor $\mathcal T$. From the definition of the E-characteristic polynomial (Definition \ref{def-2}), we see that $\mathcal T$ must have an E-eigenvector. The proof is complete. \ep

This gives an answer to Question 2 for generic tensors.

\section{E-Eigenvectors of Nonsingular Symmetric Tensors}
\setcounter{Theorem}{0} \setcounter{Proposition}{0}
\setcounter{Corollary}{0} \setcounter{Lemma}{0}
\setcounter{Definition}{0} \setcounter{Remark}{0}
\setcounter{Algorithm}{0}  \setcounter{Example}{0}
\hspace{4mm}
In this section, we consider the E-eigenvectors of symmetric tensors. A tensor $\mathcal T\in\otimes^m(\mathbb C^{n+1})$ is called symmetric if $\mathcal T\in\mathcal S^m(\mathbb C^{n+1})$, the symmetric subspace of $\otimes^m(\mathbb C^{n+1})$. By Theorem \ref{thm-8}, any nonsingular symmetric tensor must have at least one E-eigenvector. Actually, the E-eigenvectors of nonsingular symmetric tensors can be characterized by the singular points of a set of parameterized hypersurfaces.

For every given $\mathcal T\in\mathcal S^m(\mathbb C^{n+1})$, we can associated it a hypersurface with parameter $\lambda$ as
\begin{eqnarray}\label{hy}
p(\mathbf x):=\frac{1}{m}\mathbf x^T(\mathcal T\mathbf x^{m-1})-\frac{\lambda}{2}\mathbf x^T\mathbf x-(\frac{1}{m}-\frac{1}{2})\lambda.
\end{eqnarray}
In \cite{cs}, Cartwright and Sturmfels show that for nonzero $\lambda$ the corresponding E-eigenvectors are the singular points of the affine hypersurface defined by $p$. Let $\mbox{Dis}(p)$ be the classical multivariate discriminant of $p$ \cite{gkz}. It is easy to see that $\mbox{Dis}(p)$ is an univariate polynomial in the variable $\lambda$ for every fixed $\mathcal T$.

Since the singular points of the affine hypersurface is related to the discriminant of $p$. They investigate the discriminant of $p$ and show in \cite[Corollary 5.4]{cs} that $\chi_{\mathcal T}(\lambda)$ is a factor of $\mbox{Dis}(p)$. For $n=1$, they give that
$\mbox{Dis}(p)=\lambda^4\chi_{\mathcal T}(\lambda)$ when $m=3$ and $\mbox{Dis}(p)=\lambda^9\left[\chi_{\mathcal T}(\lambda)\right]^2\mbox{Det}(\mathcal T)$ when $m=4$. They point out that it would be interesting to determine the analogous factorization for arbitrary $m$ and $n$. In the following theorem, we establish the formulae for $\mbox{Dis}(p)$ for any $\mathcal T\in\mathcal S^m(\mathbb C^{n+1})$.
To this end, we establish a lemma first.

\begin{Lemma}\label{lem-00000}
Let $\mathcal T\in\otimes^m(\mathbb C^{n+1})$ and $\chi_{\mathcal T}(\lambda)$ be defined
by \reff{cha}. Then,
\begin{eqnarray}\label{dis}
\mbox{Res}\left(\left\{\begin{array}{l}\mathcal T\mathbf x^{m-1}-\lambda t^{m-2}\mathbf x\\\mathbf x^T\mathbf x-t^2\end{array}
\right.\right)
\end{eqnarray}
is equal to $\chi_{\mathcal T}(\lambda)$ when $m$ is odd and $[\chi_{\mathcal T}(\lambda)]^2$ when $m$ is even.
\end{Lemma}

\noindent {\bf Proof.} The odd case is just the definition. It is sufficient to prove the even case.

When $m$ is even, we see that $\lambda$ is a root of $\chi_{\mathcal T}(\lambda)=0$ if and only if the resultant in \reff{dis} is zero. Viewing $\lambda$ as a variable, by Lemma \ref{lem-00} which says that $\chi_{\mathcal T}(\lambda)$ is irreducible, we have that the resultant in \reff{dis} is a power of  $\chi_{\mathcal T}(\lambda)$. Comparing the degree of $\lambda$ with the help of Theorem \ref{thm-1} (iv) and \cite[Proposition 13.1.1]{gkz}, we see that the power is a square. \ep

\begin{Theorem}\label{thm-6}
Let $\mathcal T\in\mathcal S^m(\mathbb C^{n+1})$, we have
\begin{eqnarray*}
\mbox{Dis}(p)=\left\{\begin{array}{ll}\lambda^{(m-1)^{n+1}}\left[\mbox{Det}(\mathcal T)\right]^{m-3}\chi_{\mathcal T}(\lambda)&when\; m\; is\; odd,\\
\lambda^{(m-1)^{n+1}}\left[\mbox{Det}(\mathcal T)\right]^{m-3}[\chi_{\mathcal T}(\lambda)]^2&when\; m\; is\; even.\end{array}\right.
\end{eqnarray*}
\end{Theorem}

\noindent {\bf Proof.} Homogenize $p$ as
\begin{eqnarray}\label{hy-1}
q(\mathbf x,t):=\frac{1}{m}\mathbf x^T(\mathcal T\mathbf x^{m-1})-t^{m-2}\frac{\lambda}{2}\mathbf x^T\mathbf x-t^m(\frac{1}{m}-\frac{1}{2})\lambda.
\end{eqnarray}
Then, $[(\mathbf x,t)]$ is a singular point of the projective hypersurface $\{[(\mathbf x,t)]\in\mathbb P^{n+1}\;|\;q(\mathbf x,t)=0\}$ if and only if
\begin{eqnarray*}
\left\{\begin{array}{l}q(\mathbf x,t)=\frac{1}{m}\mathbf x^T(\mathcal T\mathbf x^{m-1})-t^{m-2}\frac{\lambda}{2}\mathbf x^T\mathbf x-t^m(\frac{1}{m}-\frac{1}{2})\lambda=0,\\
\mathcal T\mathbf x^{m-1}-\lambda t^{m-2}\mathbf x=0,\\
\frac{m-2}{2}t^{m-3}\lambda \left(\mathbf x^T\mathbf x-t^2\right)=0.\end{array}
\right.
\end{eqnarray*}

If $t=0$, then $\mathbf x\neq 0$ and consequently, $\mathcal T\mathbf x^{m-1}=0$. Hence, $\mbox{Det}(\mathcal T)=0$. So, if $\mbox{Det}(\mathcal T)\neq 0$, then $t\neq 0$. Let $[(\mathbf x, t)]$ be a singular point of $q$, then $\mathbf x/t$ would be a singular point of $p$ and vice verse. From \cite{gkz}, we know that the discriminant of $p$, viewed as a polynomial in the variables $\mathcal T$ and $\lambda$, has the same degree of the resultant of the derivative equations of $q$. These, together with the fact that $\mathbb C[\mathcal T,\lambda]$ is a unique factorization domain, imply that $\mbox{Dis}(p)$ is equal to the resultant of the derivative polynomial system of $q$ when $\mathcal T$ is nonsingular. Since nonsingularity is a generic hypothesis, it follows that they equal for all tensors. Actually, the discriminant of a polynomial is always equivalent to the discriminant of its homogenization, please see \cite[Chapter 13]{gkz} for more details.

Hence, we have
\begin{eqnarray*}
\mbox{Dis}(p)&=&\mbox{Res}\left(\left\{\begin{array}{l}\mathcal T\mathbf x^{m-1}-\lambda t^{m-2}\mathbf x\\
\frac{m-2}{2}t^{m-3}\lambda \left(\mathbf x^T\mathbf x-t^2\right)\end{array}
\right.\right)\\
&=&\left[\mbox{Res}\left(\left\{\begin{array}{l}\mathcal T\mathbf x^{m-1}-\lambda t^{m-2}\mathbf x\\
t\end{array}
\right.\right)\right]^{m-3}\cdot\mbox{Res}\left(\left\{\begin{array}{l}\mathcal T\mathbf x^{m-1}-\lambda t^{m-2}\mathbf x\\
\lambda(\mathbf x^T\mathbf x-t^2)\end{array}
\right.\right)\\
&=&\left[\mbox{Det}(\mathcal T)\right]^{m-3}\lambda^{(m-1)^{n+1}}\mbox{Res}\left(\left\{\begin{array}{l}\mathcal T\mathbf x^{m-1}-\lambda t^{m-2}\mathbf x\\\mathbf x^T\mathbf x-t^2\end{array}
\right.\right).
\end{eqnarray*}
Here, the second equality follows from \cite[Theorem 3.3.2(b)]{clo}, and the third from the definition of the determinant (\cite[Definition 1.2]{hhlq11}) and \cite[Proposition 13.1.1]{gkz}.
Consequently, by Lemma \ref{lem-00000}, the proof is complete.  \ep

By Lemma \ref{lem-00}, we see that the irreducible factors of $\mbox{Dis}(p)$ are $\lambda$, $\mbox{Det}(\mathcal T)$ and $\chi_{\mathcal T}(\lambda)$. Hence, Theorem \ref{thm-6} complements \cite[Corollary 5.4 and Example 5.5]{cs}. Moreover, we have the following characterization of the E-eigenvectors of a nonsingular symmetric tensor.

\begin{Corollary}\label{cor-3}
If $\mathcal T\in\mathcal S^m(\mathbb C^{n+1})$ is nonsingular, then its E-eigenvectors of an E-eigenvalue $\lambda$ are exactly the singular points of the affine hypersurface defined by $p(\mathbf x)$ with $\lambda$.
\end{Corollary}

\noindent {\bf Proof.} We have that $\lambda$ is a solution of $\mbox{Dis}(p)=0$ if and only if $p$ has a singular point $\mathbf x$ to this $\lambda$. This is further equivalent to a solution to the following system of equations:
\begin{eqnarray*}
\left\{\begin{array}{l}p(\mathbf x)=\frac{1}{m}\mathbf x^T(\mathcal T\mathbf x^{m-1})-\frac{\lambda}{2}\mathbf x^T\mathbf x-(\frac{1}{m}-\frac{1}{2})\lambda=0,\\
\mathcal T\mathbf x^{m-1}-\lambda\mathbf x=0.\end{array}
\right.
\end{eqnarray*}
Consequently, we have
\begin{eqnarray*}
0=\frac{1}{m}\mathbf x^T(\mathcal T\mathbf x^{m-1})-\frac{\lambda}{2}\mathbf x^T\mathbf x-(\frac{1}{m}-\frac{1}{2})\lambda=\frac{1}{m}\lambda\mathbf x^T\mathbf x-\frac{\lambda}{2}\mathbf x^T\mathbf x-(\frac{1}{m}-\frac{1}{2})\lambda.
\end{eqnarray*}
Since $\frac{1}{m}-\frac{1}{2}\neq 0$, we have $\mathbf x^T\mathbf x=1$. Consequently, $\mathbf x$ is an E-eigenvector with E-eigenvalue $\lambda$ if and only if it is a singular point of the affine hypersurface determined by $p$ with $\lambda$ being a root of $\mbox{Dis}(p)$.

When $\mathcal T$ is nonsingular, $0$ cannot be an E-eigenvalue of $\mathcal T$ by Theorem \ref{thm-1} (i). Hence, by Theorem \ref{thm-6}, $p$ has a singular point if and only if $\lambda\neq 0$ is a root of the E-characteristic polynomial of $\mathcal T$. Since $\mathcal T$ is nonsingular, $\lambda$ is an E-eigenvalue of $\mathcal T$ if and only if it is a root of the E-characteristic polynomial by Theorem \ref{thm-1} (ii). Combining the two results, the corollary follows. \ep

\section{Final Remarks}
\setcounter{Theorem}{0} \setcounter{Proposition}{0}
\setcounter{Corollary}{0} \setcounter{Lemma}{0}
\setcounter{Definition}{0} \setcounter{Remark}{0}
\setcounter{Algorithm}{0}  \setcounter{Example}{0}
\hspace{4mm} In \cite{hhlq11}, the determinant of a tensor is demonstrated to be an important function in the eigenvalue theory of tensors. In \reff{ideal}, a certain power of the determinant of a tensor is expressed as a polynomial combination of orthogonal invariants of the underlying tensor. This sheds light on further investigation for the determinant using invariant theory.

Actually, we have some preliminary observations from another perspective here.
For a third order ($n+1$)-dimensional tensor $\mathcal T$, we associated it a set of $n+1$ matrices $\{A^{(0)},\ldots,A^{(n)}\}$ such that
\begin{eqnarray*}
(\mathcal T\mathbf x^2)_i=\mathbf x^TA^{(i)}\mathbf x,\;\forall i\in(n),
\end{eqnarray*}
where $a^{(i)}_{jk}=t_{ijk},\;\forall j,k\in (n)$.

Let $\mathbb O(n+1)$ be the $(n+1)\times (n+1)$ orthogonal matrices with entries in $\mathbb C$. Define the group action of $\mathbb O(n+1)$ on the set $\otimes^3(\mathbb C^{n+1})$ of third order ($n+1$)-dimensional tensors as:
\begin{eqnarray*}
(G\cdot\mathcal T)_{ijk}:=\sum_{p,q\in(n)}t_{ipq}g_{jp}g_{kq},\;\forall i,j,k\in(n).
\end{eqnarray*}
It is equivalent to the simultaneous conjugation action on the set of matrices $\{A^{(0)},\ldots,A^{(n)}\}$, i.e., $GA^{(i)}G^T$. It is different from the orthogonal group action defined in \cite{q,q2,lqz}, which was used in the previous sections. By
\cite[Section 3.1.4]{dm}, we see that $\mbox{Det}(\mathcal T)$ is an invariant under this action as well. The invariants of this group action were investigated extensively in the literature.
We have the following result, which is due to Sibirskii \cite{s68} and Procesi \cite{p76}.
\begin{Theorem}\label{thm-9}
Let $\mathcal T\in \otimes^3(\mathbb C^{n+1})$ and matrices $\{A^{(0)},\ldots,A^{(n)}\}$ be defined as above. Then, every invariant under the group action defined above can be expressed as a combination of a finite number of polynomials
\begin{eqnarray*}
\mbox{Tr}\left[p_j(A^{(0)},\ldots,A^{(n)},(A^{(0)})^T,\ldots,(A^{(n)})^T)\right], j\in[N]
\end{eqnarray*}
with polynomial coefficients for some positive integer $N$.
\end{Theorem}

We can assume that all of the matrices $\{A^{(0)},\ldots,A^{(n)}\}$ are symmetric, without loss of generality. When $n=1$, we have the following result, whose proof follows from \cite[Corollary 7.8]{hhlq11} by direct computation.

\begin{Theorem}\label{thm-10}
Let $\mathcal T\in \otimes^3(\mathbb C^{2})$ and $\mathcal T=[A,B]$ with $A$ and $B$ being symmetric. Then, we have
\begin{eqnarray*}
\mbox{Det}(\mathcal T)&=&\left[\mbox{Tr}(A)\mbox{Tr}(B)-\mbox{Tr}(AB)+\mbox{Tr}(A)^2-\mbox{Tr}(A^2)\right]\left[\mbox{Tr}(B^2)-\mbox{Tr}(B)^2\right].
\end{eqnarray*}
\end{Theorem}

{\bf Acknowledgement.} We are grateful to Professor Mattias Jonsson and Professor Giorgio Ottaviani for
valuable suggestions.


\begin{thebibliography}{abc99xyz}

\bibitem{cs}
D. Cartwright, and B. Sturmfels, The number of eigenvalues of a
tensor, Linear Algebra Appl. 438 (2013), pp. 942--952.

\bibitem{c49}
W.-L. Chow, On compact complex analytic varieties, Amer. J. Math. 71 (1949), pp. 893--914.

\bibitem{clo1}
D. Cox, J. Little, and D. O\'{}Shea, Ideals, Varieties, and Algorithms: An
Introduction to Computational Algebraic Geometry and Commutative
Algebra, New York: Springer-Verlag, 2006.

\bibitem{clo}
D. Cox, J. Little, and D. O\'{}Shea, Using Algebraic Geometry, New
York: Springer-Verlag, 1998.

\bibitem{dm}
V. Dolotin, and A. Morozov, Introduction to Non-Linear Algebra, World Scientific, 2007.

\bibitem{fs}
J.E. Fornaess, and N. Sibony, Complex dynamics in higher dimension, I, Ast\'erisque 5 (1994), pp. 201--231.

\bibitem{fo12}
S. Friedland, and G. Ottaviani,
The number of singular vector tuples and uniqueness of best rank one approximation of tensors, arXiv:1210.8316.

\bibitem{f84}
W. Fulton, Intersection Theory, Springer, Berlin, 1984.

\bibitem{gkz}
I.M. Gelfand, M.M. Kapranov, and A.V. Zelevinsky, Discriminants, Resultants and Multidimensional Determinants, Birkh\"{a}user, Boston, 1994.

\bibitem{gh}
P. Griffiths, and J. Harris, Priniciples of Algebraic Geometry, Wiley, 1978.

\bibitem{h77}
R. Hartshorne, Algebraic Geometry, Graduate Texts in Mathematics 52,
Springer, New York, 1977.

\bibitem{h12}
S. Hu, E-characteristic polynomials tensors, Talk at the RTG workshop on Tensors and their Geometry in High Dimensions, September 26--September 29, 2012, MSRI, University of California at Berkeley.

\bibitem{hhlq11}
S. Hu, Z.-H. Huang, C. Ling, and L. Qi, On determinants and eigenvalue theory of tensors, J. Symbolic Comput. 50 (2013), pp. 508--531.

\bibitem{hq12}
S. Hu, and L. Qi, Algebraic connectivity of an even uniform hypergraph, J. Comb. Optim. 24 (2012), pp. 564--579.

\bibitem{hq13}
S. Hu, and L. Qi, E-characteristic polynomial of a tensor of dimension two, Appl. Math. Lett. 26 (2013), pp. 225--231.

\bibitem{lqz}
A.-M. Li, L. Qi, and B. Zhang, E-characteristic polynomials of tensors, Commun. Math. Sci. 11 (2013) 33--53.

\bibitem{l05}
L.-H. Lim, Singular values and eigenvalues of tensors: a variational approach,
Computational Advances in Multi-Sensor Adaptive Processing, 2005 1st IEEE
International Workshop on, 2005, pp. 129--132.

\bibitem{nqww}
G. Ni, L. Qi, F. Wang, and Y. Wang, The degree of the E-characteristic polynomial of
an even order tensor, J. Math. Anal. Appl., 329 (2007), pp. 1218--1229.

\bibitem{oo12}
L. Oeding, G. Ottaviani, Eigenvectors of tensors and algorithms for Waring
decomposition, J. Symbolic Comput., to appear.

\bibitem{oss80}
C. Okonek, M. Schneider, H. Spindler, Vector bundles on complex projective spaces, Progress in Mathematics, vol. 3,
Birkh\"{a}user Boston, Mass., 1980.

\bibitem{os12}
G. Ottaviani, B. Sturmfels, Matrices with Eigenvectors in a Given Subspace,
Proc. of the American Math. Soc., in press.

\bibitem{p76}
C. Procesi, The invariant theory of $n\times n$ matrices, Adv. Math., 19 (1976), pp. 306-381.

\bibitem{q}
L. Qi, Eigenvalues of a real supersymmetric tensor, J. Symbolic Comput. 40 (2005), pp. 1302--1324.

\bibitem{q2}
L. Qi, Eigenvalues and invariants of tensors, J. Math. Anal. Appl. 325 (2007), pp. 1363--1377.

\bibitem{s68}
K. S. Sibirskii, Algebraic invariants for a set of matrices, Siberian Math. J. 9 (1968), pp. 115-124.

\bibitem{skkt00}
K.E. Smith, L. Kahanp\"{a}\"{a}, P. Kek\"{a}l\"{a}inen, W. Traves, An Invitation to Algebraic Geometry, Springer, New York, 2000.

\end{thebibliography}
\end{document}